\input amstex
\documentstyle{amsppt}
\pagewidth{6.0in}\vsize8.5in\parindent=6mm
\parskip=3pt\baselineskip=14pt\tolerance=10000\hbadness=500
\document
\topmatter
\title $\zeta(5)$ is irrational
\endtitle
\author Yong-Cheol Kim \endauthor
\abstract We present an elementary proof of the irrationality of
$\zeta(5)$ based upon the Dirichlet's approximation theorem and the
Prime Number Theorem.
\endabstract
\thanks AMS 2000 Mathematics Subject Classification. 11J81, 11M06.
\endthanks
\address Department of Mathematics Education, Korea University, Seoul 136-701, Korea
\endaddress
\email ychkim$\@$korea.ac.kr
\endemail
\endtopmatter

\define\inn#1#2{\langle#1,#2\rangle}

\define\lcontr{\rfloor}
\define\lco#1#2{{#1}\lcontr{#2}}
\define\lcoi#1#2{\imath({#1}){#2}}
\define\rco#1#2{{#1}\rcontr{#2}}

\define\ap{\alpha}             
\define\bt{\beta}
\define\gm{\gamma}             
\define\dt{\delta}             
\define\vep{\varepsilon}
\define\zt{\zeta}

            \define\iy{\infty}
\define\lt{\left}            \define\rt{\right}
\define\f{\frac}             \define\el{\ell}


\define\fm{{\frak m}}

\define\BN{{\Bbb N}}

\define\BQ{{\Bbb Q}}
\define\BR{{\Bbb R}}

\define\BZ{{\Bbb Z}}

\define\cI{{\Cal I}}




\define\n{\nabla}            \define\e{\eta}
\define\pa{\partial}        \define\fd{\fallingdotseq}
       
     \define\ds{\dsize} 
\define\pf{\noindent{\it Proof. }}  \define\rM{\text{\rm M}}
\define\rP{\text{\rm P}}      \define\rQ{\text{\rm Q}}

\subheading{1. Introduction}

We are very concerned about a question of an arithmetic nature of
the values of the Riemann zeta function
$$\zt(z)=\sum_{n=1}^{\iy}\f{1}{n^z}$$ at integral points
$z=2,3,4,5,\cdots\,$, which has been a challenge in Number Theory
area. Originally L. Euler obtained the exact values of $\zt(z)$ at
all even integral numbers $z=2,4,6,\cdots$ as follows;
$$\zt(2k)=(-1)^{k-1} 2^{2k-1} c_{2k-1}\f{\pi^{2k}}{(2k-1)!},\,k\in\BN,$$
where the sequence $c_{\el}$ of rational numbers is defined by the
Laurent expansion $$\f{1}{e^z
-1}=\f{1}{z}-\f{1}{2}+\sum_{\el=1}^{\iy}\f{c_{\el}}{\el!}\,z^k.$$ He
also made very serious attempts to evaluate $\zt(2k+1),k\in\BN,$ and
calculated to 15 places of decimals along with the corresponding
quotients $\zt(2k+1)\pi^{-(2k+1)}$. All that is known about these
values is surprisingly recent result of Ap\'ery (1978) ( see [1] and
[3] ) that $\zt(3)$ is irrational. Along with this problem, an
interesting question about the transcendentality of
$\zt(2k+1),k\in\BN$ is still far from being solved. However some
progress have been made on the irrationality of $\zt(2k+1)$. In
2000, K. Ball and T. Rivoal ( see [2] and [4] ) proved that
infinitely many of the numbers $\zt(3),\zt(5),\zt(7),\cdots$ are
irrational. Also W. Zudilin proved an interesting result in 2001 (
see [5] and [6] ) that at least one of the four numbers
$\zt(5),\zt(7),\zt(9),\zt(11)$ is irrational.

In this paper, we shall furnish an elemetary proof of the
irrationality of $\zt(5)$ based upon the Dirchlet's approximation
theorem and the Prime Number Theorem.

\proclaim{Theorem 1.1} $\zt(5)$ is irrational. \endproclaim

\subheading{2. Preliminary estimates}

First of all, we recall the Dirichlet's approximation theorem and
the Prime Number Theorem to be used as important tools for our
proof.

\proclaim{[Dirichlet's Approximation Theorem]} \noindent For any
$\ap\in\BR$ and $N\in\BN$, there are $n\in\BN\,\,( n\le N )$ and
$p\in\BZ$ such that $\,\,\ds\bigl|\ap-\f{p}{n}\bigr|<\f{1}{N
n}.$\endproclaim

We can easily deduce the following key fact which is a variant of
the Dirichlet's approximation theorem.

\proclaim{[Key Lemma]} $\ap\in\BQ^c$ if and only if for any
$\vep>0$, there are some $x\in\BN$ and $y\in\BZ$ such that $0<|\ap
x-y|<\vep$.\endproclaim


\proclaim{[Prime Number Theorem]} If $\pi(n)$ denotes the number of
primes $p\le n$, then $\,\ds\lim_{n\to\iy}\f{\pi(n)}{n/\ln
n}=1.$\endproclaim

Next we shall give various useful technical lemmas to be shown by
using partial fractions without detailed proof. We observe that
for $r=(r_1,r_2,r_3,r_4,r_5)\in (\BN\cup\{0\})^5$,
$$ \split \cI(r)&\fd\idotsint_{[0,1]^5}
\f{x_1^{r_1}x_2^{r_2}x_3^{r_3}x_4^{r_4}x_5^{r_5}}{1-x_1 x_2 x_3
x_4 x_5}\,dx_1\, dx_2\, dx_3\, dx_4\, dx_5
\\&=\sum_{k=0}^{\iy}
\f{1}{(k+r_1+1)(k+r_2+1)(k+r_3+1)(k+r_4+1)(k+r_5+1)}\endsplit
\tag{2.1}$$ by expanding $(1-x_1 x_2 x_3 x_4 x_5)^{-1}$ as a
geometric series.

\proclaim{Lemma 2.1} If $\,r_1>r_2>r_3>r_4>r_5\ge 0$, then we have
that
$$\split
\cI(r)&=\f{1}{(r_1-r_5)(r_2-r_5)(r_3-r_5)(r_4-r_5)}\lt(\f{1}{r_5
+1}+\f{1}{r_5 +2}+\cdots +\f{1}{r_4}\rt) \\
&-\f{1}{(r_1-r_5)(r_3-r_4)}\lt(\f{1}{(r_2-r_5)(r_3-r_5)}+\f{1}{(r_2-r_5)(r_2-r_4)}
+\f{1}{(r_1-r_4)(r_2-r_4)}\rt) \\
&\qquad\qquad\qquad\qquad\times \lt(\f{1}{r_4 +1}+\f{1}{r_4
+2}+\cdots+\f{1}{r_3}\rt) \\
&+\f{1}{(r_1-r_5)(r_2-r_3)}\lt(\f{1}{(r_2-r_5)(r_2-r_4)}+\f{1}{(r_1-r_4)(r_2-r_4)}
+\f{1}{(r_1-r_4)(r_1-r_3)}\rt) \\
&\qquad\qquad\qquad\qquad\times\lt(\f{1}{r_3 +1}+\f{1}{r_3
+2}+\cdots+\f{1}{r_2}\rt) \\
&-\f{1}{(r_1-r_5)(r_1-r_4)(r_1-r_3)(r_1-r_2)}\lt(\f{1}{r_2
+1}+\f{1}{r_2 +2}+\cdots+\f{1}{r_1}\rt). \endsplit $$
\endproclaim

\proclaim{Lemma 2.2} If $\,r_1=r_2>r_3>r_4>r_5\ge 0$, then we have
that
$$\split \cI(r)&=\f{1}{(r_2-r_5)^2
(r_3-r_5)(r_4-r_5)}\lt(\f{1}{r_5 +1}+\f{1}{r_5 +2}+\cdots+\f{1}{r_2}\rt) \\
&-\f{1}{(r_2-r_5)(r_2-r_4)}\lt(\f{1}{(r_3-r_4)(r_3-r_5)}+\f{1}{(r_3-r_4)(r_2-r_4)}+\f{1}{(r_3-r_5)(r_4-r_5)}\rt)
\\ &\qquad\qquad\qquad\qquad\times\lt(\f{1}{r_4 +1}+\f{1}{r_4 +2}+\cdots+\f{1}{r_2}\rt) \\
&+\f{1}{(r_2-r_5)(r_2-r_3)}\lt(\f{1}{(r_3-r_4)(r_3-r_5)}+\f{1}{(r_3-r_4)(r_2-r_4)}+\f{1}{(r_2-r_3)(r_2-r_4)}\rt)
\\ &\qquad\qquad\qquad\qquad\times\lt(\f{1}{r_3 +1}+\f{1}{r_3 +2}+\cdots+\f{1}{r_2}\rt) \\
&-\f{1}{(r_2-r_5)(r_2-r_4)(r_2-r_3)}\lt(\zt(2)-\lt(1+\f{1}{2^2}+\f{1}{3^2}+\cdots+\f{1}{r_2^2}\rt)\rt).
\endsplit $$ \endproclaim

\proclaim{Lemma 2.3} If $\,r_1=r_2=r_3>r_4>r_5\ge 0$, then we have
that
$$\split \cI(r)&=\f{1}{(r_3-r_5)^3 (r_4-r_5)}\lt(\f{1}{r_5
+1}+\f{1}{r_5 +2}+\cdots+\f{1}{r_3}\rt) \\
&-\f{1}{(r_3-r_5)(r_3-r_4)^2}\lt(\f{1}{r_4-r_5}+\f{1}{r_3-r_4}\rt)\lt(\f{1}{r_4
+1}+\f{1}{r_4 +2}+\cdots+\f{1}{r_3}\rt) \\
&+\f{1}{r_3-r_5}\lt(\f{1}{(r_3-r_4)(r_4-r_5)}+\f{1}{(r_3-r_4)^2}
-\f{1}{(r_3-r_5)(r_4-r_5)}\rt) \\
&\qquad\qquad\qquad\qquad\times\lt(\zt(2)-\lt(1+\f{1}{2^2}+\f{1}{3^2}+\cdots+\f{1}{r_3^2}\rt)\rt)
\\
&+\f{1}{(r_3-r_5)(r_3-r_4)}\lt(\zt(3)-\lt(1+\f{1}{2^3}+\f{1}{3^3}+\cdots+\f{1}{r_3^3}\rt)\rt).
\endsplit $$ \endproclaim

\proclaim{Lemma 2.4} If $\,r_1=r_2=r_3=r_4>r_5\ge 0$, then we have
that
$$\split \cI(r)&=\f{1}{(r_4-r_5)^4}\lt(\f{1}{r_5 +1}+\f{1}{r_5
+2}+\cdots+\f{1}{r_4}\rt) \\
&-\f{1}{(r_4-r_5)^3}\lt(\zt(2)-\lt(1+\f{1}{2^2}+\f{1}{3^2}+\cdots+\f{1}{r_4^2}\rt)\rt)
\\
&-\f{1}{(r_4-r_5)^2}\lt(\zt(3)-\lt(1+\f{1}{2^3}+\f{1}{3^3}+\cdots+\f{1}{r_4^3}\rt)\rt)
\\
&-\f{1}{r_4-r_5}\lt(\zt(4)-\lt(1+\f{1}{2^4}+\f{1}{3^4}+\cdots+\f{1}{r_4^4}\rt)\rt).
\endsplit $$ \endproclaim

\proclaim{Lemma 2.5} If $\,r_1=r_2>r_3>r_4=r_5\ge 0$, then we have
that
$$\split \cI(r)&=\f{1}{(r_2-r_4)^2
(r_3-r_4)}\lt(\zt(2)-\lt(1+\f{1}{2^2}+\f{1}{3^2}+\cdots+\f{1}{r_4^2}\rt)\rt)
\\ &-\f{1}{(r_2-r_4)^2
(r_3-r_4)}\lt(\f{1}{r_3-r_4}+\f{1}{r_2-r_3}\rt)\lt(\f{1}{r_4
+1}+\f{1}{r_4 +2}+\cdots+\f{1}{r_3}\rt) \\
&+\f{1}{(r_2-r_4)^3}\lt(\f{1}{r_2-r_3}-\f{1}{r_3-r_4}\rt)\lt(\f{1}{r_4
+1}+\f{1}{r_4 +2}+\cdots+\f{1}{r_3}\rt) \\
&\f{1}{(r_2-r_4)^2
(r_2-r_3)}\lt(\f{1}{r_2-r_3}+\f{1}{r_3-r_4}\rt)\lt(\f{1}{r_3
+1}+\f{1}{r_3 +2}+\cdots+\f{1}{r_2}\rt) \\
&-\f{1}{(r_2-r_4)^2
(r_2-r_3)}\lt(\zt(2)-\lt(1+\f{1}{2^2}+\f{1}{3^2}+\cdots+\f{1}{r_2^2}\rt)\rt).
\endsplit $$ Here we note that
$1+1/2^2+1/3^2+\cdots+1/r_4^2$ could be regarded as $0$ if
$\,r_4=0$.
\endproclaim

\proclaim{Lemma 2.6} If $\,r_1=r_2=r_3>r_4=r_5\ge 0$, then we have
that
$$\split
\cI(r)&=\f{1}{(r_3-r_4)^3}\lt(\zt(2)-\lt(1+\f{1}{2^2}+\f{1}{3^2}+\cdots+\f{1}{r_4^2}\rt)\rt)
\\ &-\f{3}{(r_3-r_4)^4}\lt(\f{1}{r_4 +1}+\f{1}{r_4
+2}+\cdots+\f{1}{r_3}\rt) \\
&+\f{2}{(r_3-r_4)^3}\lt(\zt(2)-\lt(1+\f{1}{2^2}+\f{1}{3^2}+\cdots+\f{1}{r_3^2}\rt)\rt)
\\
&+\f{1}{(r_3-r_4)^2}\lt(\zt(3)-\lt(1+\f{1}{2^3}+\f{1}{3^3}+\cdots+\f{1}{r_3^3}\rt)\rt).
\endsplit $$ Here we note that
$1+1/2^2+1/3^2+\cdots+1/r_4^2$ could be regarded as $0$ if
$\,r_4=0$. \endproclaim

\proclaim{Lemma 2.7} If $\,r_1=r_2=r_3=r_4=r_5\ge 0$, then we have
that
$$\cI(r)=\zt(5)-\lt(1+\f{1}{2^5}+\f{1}{3^5}+\cdots+\f{1}{r_1^5}\rt).$$
Here we note that $1+1/2^5+1/3^5+\cdots+1/r_1^5$ could be regarded
as $0$ if $\,r_1=0$.\endproclaim

\subheading{3. The proof of Theorem 1.1}

In this section, we shall prove Theorem 1.1. From Lemmas 2.1, 2.2,
2.3, 2.4, 2.5, 2.6, and 2.7, we can easily obtain the following
proposition.

\proclaim{Proposition 3.1} If we denote by $\fm(n)$ the least
common multiple of $1,2,\cdots,n$, then we have that
$$\split \cI_n&\fd\idotsint_{[0,1]^5}
\f{(1-x_1)^n (1-x_2)^n (1-x_3)^n (1-x_4)^n\,\rP_n(x_5)}{1-x_1 x_2
x_3
x_4 x_5}\,dx_1\, dx_2\, dx_3\, dx_4\,dx_5 \\
&=\f{a_n\zt(2)+b_n\zt(3)+c_n\zt(4)+d_n\zt(5)+e_n}{{\fm(n)}^5}
\endsplit$$ where $a_n,b_n,c_n,d_n,e_n$ are some integers and
$\rP_n(x_5)=\ds\f{1}{n!}\lt(\f{d}{dx_5}\rt)^n (x_5^n (1-x_5)^n)$.
\endproclaim

Applying the integration by parts $n$-times with respect to $x_5$ to
the integral $\cI_n$ leads us to get
$$\cI_n=\idotsint_{[0,1]^5}
\f{x_1^n(1-x_1)^n x_2^n(1-x_2)^n x_3^n(1-x_3)^n x_4^n(1-x_4)^n
x_5^n(1-x_5)^n}{(1-x_1 x_2 x_3 x_4 x_5)^{n+1}}\,dx_1\, dx_2\, dx_3\,
dx_4\,dx_5.$$ We now consider the function $\rQ$ on $[0,1]^5$
defined by $$\rQ(x)=\f{x_1(1-x_1) x_2(1-x_2) x_3(1-x_3) x_4(1-x_4)
x_5(1-x_5)}{1-x_1 x_2 x_3 x_4 x_5}$$ where $x$ means a multiindex
$x=(x_1,x_2,x_3,x_4,x_5)\in [0,1]^5$. Then we shall try to obtain
the nice upper bound of the function $\rQ(x)$ on $[0,1]^5$ which is
not its maximum value, but suitable for our goal.

\proclaim{Lemma 3.2} The function $\rQ(x)$ satisfies the following
property;
$$\sup_{x\in [0,1]^5}
|\rQ(x)|\le\f{\ds\lt(\f{1}{2}+\f{1}{100}\rt)^5
\lt(\f{1}{2}-\f{1}{200}\rt)^5}{1-\ds\lt(\f{1}{2}+\f{1}{100}\rt)^5}
=\f{(0.25245)^5}{0.9654974749}\fd\gm\,.$$ \endproclaim

\noindent{\it Proof.} Since the function $\rQ(x)$ is continuous on
$[0,1]^5$, differentiable infinitely on $(0,1)^5$, and vanishes on
the boundary $\pa([0,1]^5)$ of $[0,1]^5$, it has the maximum value
at some point $x^0=(x_1^0,x_2^0,x_3^0,x_4^0,x_5^0)\in (0,1)^5$, and
so $x^0$ is a critical point for $\rQ$. It follows from simple
calculation that $\n\rQ(x)=0$ on $(0,1)^5$ if and only if
$$\split 1+x_1(x_1 x_2 x_3 x_4 x_5-2)&=0, \\
         1+x_2(x_1 x_2 x_3 x_4 x_5-2)&=0, \\
         1+x_3(x_1 x_2 x_3 x_4 x_5-2)&=0, \\
         1+x_4(x_1 x_2 x_3 x_4 x_5-2)&=0, \\
         1+x_5(x_1 x_2 x_3 x_4 x_5-2)&=0\,\,\quad\text{on $(0,1)^5$.}\endsplit $$
This implies that if $x\in (0,1)^5$ is an critical point for $\rQ$,
then it should satisfy the condition $x_1=x_2=x_3=x_4=x_5$. Thus in
order to trace out the critical points it is natural for us to
consider the function $F(t)$ on $(0,1)$ defined by $$F(t)\fd
\rQ(t,t,t,t,t)=\f{t^5 (1-t)^5}{1-t^5}.\tag{3.1}$$ Thus we shall
track down the critical points for $F(t)$ instead of doing those for
$\rQ(x)$. We observe that $F'(t)=0$ on $(0,1)$ if and only if
$1-2t+t^6=0$. If we set $G(t)=1-2t+t^6$ on $(0,1)$, then we have
that $$G'(t_0)=-2 + 6
t_0^5=0\,\,\,\Leftrightarrow\,\,\,t_0=\lt(\f{1}{3}\rt)^{1/5}.$$ Thus
we easily see that $G(t)$ is decreasing on $(0,t_0]$, increasing on
$[t_0,1)$, $\lim_{t\to 0^+} G(t)=1$, $G(t_0)=-0.33790260\cdots<0$,
and $\lim_{t\to 1^-}G(t)=0$. Hence we see that there exist only one
critical point $t_1\in (0,1)$ for $F$ ( i.e. $F'(t_1)=0$ ) and also
we can expect that the point $t_1$ is near $t=1/2$ because
$G(1/2)=1/2^6$. In fact, it follows from simple computation that
$$G\lt(\f{1}{2}+\f{1}{200}\rt)=0.006586252353140625>0\,\,\,\,\text{ and
}\,\,\,\,G\lt(\f{1}{2}+\f{1}{100}\rt)=-0.002403712199<0.$$ This
implies that $$\f{1}{2}+\f{1}{200}\le
t_1\le\f{1}{2}+\f{1}{100}\tag{3.2}$$ Therefore by (3.1) and (3.2) we
can conclude that
$$\sup_{x\in [0,1]^5}
|\rQ(x)|=\sup_{t\in (0,1)}F(t)=\sup_{t\in
[\f{1}{2}+\f{1}{200},\f{1}{2}+\f{1}{100}]}F(t)\le\f{\ds\lt(\f{1}{2}+\f{1}{100}\rt)^5\cdot
\lt(\f{1}{2}-\f{1}{200}\rt)^5}{1-\ds\lt(\f{1}{2}+\f{1}{100}\rt)^5}.\,\,\,\,\qed$$

From Proposition 3.1 and Lemma 3.2, we can easily obtain the
following lemma.

\proclaim{Lemma 3.3} If $\,a_n,b_n,c_n,d_n,e_n,$ and $\gm$ are the
integers given in Proposition 3.1 and Lemma 3.2, then we have that
$$0<\e_n\le a_n\zt(2)+b_n\zt(3)+c_n\zt(4)+d_n\zt(5)+e_n\le
(3^5\gm)^n\cdot\zt(5)\fd\dt^n\cdot\zt(5)$$ for all sufficiently
large $n$. Here in fact it turns out that
$\dt=0.2580667226431440537\cdots<1$ and
$$\e_n=\f{\fm(n)^5 4^{-10
n}}{2^5(1-4^{-5})^{n+1}}>0.$$
\endproclaim
\pf By the Prime Number Theorem, we can derive that $\pi(n)\le\ds\ln
3\cdot\f{n}{\ln n}$. Thus we obtain that
$$\fm(n)\le n^{\pi(n)}\le n^{n(\f{\ln 3}{\ln n})}=3^n$$ for all
sufficiently large $n$. Since we see that
$$\idotsint_{[0,1]^5} \f{1}{1-x_1 x_2 x_3 x_4 x_5}\,dx_1\, dx_2\,
dx_3\, dx_4\,dx_5=\zt(5),$$ the second inequality easily follows
from Proposition 3.1 and Lemma 3.2.

The first inequality can be obtained from Proposition 3.1 and the
following inequality
$$\split\f{a_n\zt(2)+b_n\zt(3)+c_n\zt(4)+d_n\zt(5)+e_n}{\fm(n)^5}
&\ge\idotsint_{[\f{1}{4},\f{3}{4}]^5} \f{\prod_{i=1}^5
x_i^n(1-x_i)^n}{(1-x_1 x_2 x_3 x_4 x_5)^{n+1}}\,dx_1\, dx_2\, dx_3\,
dx_4\,dx_5 \\
&\ge\f{4^{-5n}4^{-5n}}{(1-4^{-5})^{n+1}}\cdot\biggl(\f{1}{2}\biggr)^5.
\endsplit$$ Hence we complete the proof. \qed

We finally prove the irrationality of $\zt(5)$ by applying Lemma 3.3
and the Dirichlet's approximation theorem.

{\bf [ Proof of Theorem 1.1. ]} We fix any $\vep>0$. For $n\in\BN$,
we set $\ap_n=a_n\zt(2)+b_n\zt(3)+c_n\zt(4)$ where
$\,a_n,b_n,c_n,d_n$ and $e_n$ are the integers given in Proposition
3.1

If there exists some sufficiently large $N_0\in\BN$ so that
$\ap_n=0$ and $\dt^n\cdot\zt(5)<\vep$ for any $n\ge N_0$, then it
follows from Lemma 3.3 that
$$0<|d_n\zt(5)+e_n|\le\dt^n\cdot\zt(5)<\vep$$ for any $n\ge N_0$. Thus
we can complete the proof in this case.

If such $N_0$ never exists ( i.e. if there exists an increasing
subsequence $\{n_k\}\subset\BN$ with $\lim_{k\to\iy} n_k=\iy$ such
that $\ap_{n_k}\neq 0$ for all $k\in\BN$ ), then by Lemma 3.3 we may
choose some sufficiently large $N_1\in\BN$ so that
$$0<\e_{n_k}\le\bt_{n_k}=\ap_{n_k}+d_{n_k}\zt(5)+e_{n_k}\le\dt^{n_k}\zt(5)$$
whenever $n_k\ge N_1$. To avoid the complexity of notations, from
now on we use $k$ and $k\ge N_1$ in place of $n_k$ and $n_k\ge N_1$,
respectively. With the notations changed, we have that
$$0<\e_k\le\bt_k=\ap_k+d_k\zt(5)+e_k\le\dt^k\zt(5)\tag{3.3}$$
whenever $k\ge N_1$. Since $\lim_{k\to\iy}\bt_k=0$ by (3.3), for
each $k\ge N_1$ we may choose a sufficiently large $L_k>0$ such that
$$\f{1}{L_k^{\tau}}<\bt_k^{1/2}<\f{1+\nu}{L_k^{\tau}}\tag{3.4}$$
where $\tau\in(1/2,1)$ and $\nu\in(0,1/2)$ are certain constants to
be determined later. For each $k\ge N_1$, we now set
$$N_k=\f{L_k^{1-\tau}}{\bt_k^{1/2}}.\tag{3.5}$$ By the
Dirichlet's approximation theorem, for each $k\ge N_1$ there exist
some $p_k\in\BZ$ and $m_k\in\BN$ with $m_k\le N_k$ such that
$$\biggl|\ap_k-\f{p_k}{m_k}\biggr|<\f{1}{N_k m_k}.\tag{3.6}$$
By (3.4), for each $k\ge N_1$ there exists some
$\gm_k\in(0,1)$ such that
$$0<\f{1}{L_k}<\gm_k<\f{1}{N_k}=\f{\bt_k^{1/2}}{L_k^{1-\tau}}\le\f{m_k}{N_k}\le
1.\tag{3.7}$$ Since $\lim_{k\to\iy}\bt_k=0$ and $1/2<\tau<1$, we
choose some sufficiently large $k_0\ge N_1$ such that
$$\max\bigl\{\bt_{k_0}^{\f{1}{2}},2\,\bt_{k_0}^{1-\f{1}{2\tau}}\bigr\}<\f{\vep}{2}.\tag{3.8}$$
Then we may select some very small $\mu\in(0,1/100)$ so that
$$\f{1}{L_{k_0}}<\f{1}{L_{k_0}^{1-\mu}}=\biggl(\f{1}{L_{k_0}}\biggr)^{1-\mu}<\gm_{k_0}<1.\tag{3.9}$$
We now choose $\tau=\f{1}{2}\bigl(1+\f{\mu}{2}\bigr)$. Since
$0<\bt_{k_0}<1$, if $\nu\in(0,1/2)$ could be chosen sufficiently
close to $0$ in the last inequality of (c) below, then by (3.4),
(3.5), (3.7), (3.8) and (3.9) we obtain three facts as follows.
$$\split (a)\,\,\,\,
&\f{1}{N_{k_0}}\le\bt_{k_0}^{\f{1}{2}}<\f{\vep}{2}.\\
(b)\,\,\,\,&m_{k_0}\bt_{k_0}\le N_{k_0}\bt_{k_0}\le
L_{k_0}^{1-\tau}\bt_{k_0}^{\f{1}{2}}<(1+\nu)^{\f{2-\mu}{2+\mu}}\bt_{k_0}^{1-\f{1}{2\tau}}\le
2\,\bt_{k_0}^{1-\f{1}{2\tau}}<\f{\vep}{2}.\\
(c)\,\,\,\,&\f{1}{N_{k_0}m_{k_0}}\le\f{1}{\gm_{k_0}N_{k_0}^2}=\f{\bt_{k_0}}{\gm_{k_0}L_k^{2(1-\tau)}}
<\f{1}{\gm_{k_0}}\bt_{k_0}^{\f{1}{\tau}}<L_{k_0}^{1-\mu}\bt_{k_0}^{\f{1}{\tau}}
=L_k^{-\mu}(1+\nu)^{\f{1}{\tau}}\bt_{k_0}^{\f{1}{2\tau}}\\
&\qquad\qquad<(1+\nu)^{\f{4}{2+\mu}}\bt_{k_0}^{\f{1+\mu}{1+\f{\mu}{2}}}<\bt_{k_0}.
\endsplit$$ We set ${\rM}_0=\ds\f{p_{k_0}}{m_{k_0}}+d_{n_{k_0}}\zt(5)+e_{n_{k_0}}$. Then by
(3.3), (3.6), (c) and the triangle inequality, we obtain that
$$\split {\rM}_0&\ge\beta_{k_0}-\frac{1}{N_{k_0}
m_{k_0}}\ge\beta_{n_{k_0}}-(1+\nu)^{\frac{4}{2+\mu}}\beta_{k_0}^{\frac{1+\mu}{1+\frac{\mu}{2}}}>0.\endsplit$$
Also it follows from (3.3), (3.6), (a), (b) and the triangle
inequality that $$\split
{\rM}_0\le\biggl|\f{p_{k_0}}{m_{k_0}}-\alpha_{k_0}\biggr|+|\alpha_{k_0}+d_{k_0}\zeta(5)+e_{k_0}|<\frac{1}{N_{k_0}
m_{k_0}}+\beta_{k_0}<\frac{\varepsilon}{2
m_{k_0}}+\frac{\varepsilon}{2 m_{k_0}}<\frac{\varepsilon}{m_{k_0}}.
\endsplit$$ Hence this implies that
$$0<m_{k_0}\,{\rM}_0=|m_{k_0}\,d_{k_0}\zeta(5)+p_{k_0}+m_{k_0}\,e_{k_0}|<\varepsilon.$$
Therefore we can complete the proof by applying the Key Lemma. \qed

{\bf Acknowledgement.} The author would like to thank Victor
Scharaschkin and Sander Zwegers for their helpful comments.

{\bf Added in the proof.} I obtained an elementary proof of the
irrationality of $\zt(2n+1)$ after I had submitted this paper
somewhere else.

\enddocument